\documentclass[12pt]{amsart}
\usepackage{amsmath,amscd,amsthm,amsfonts, amssymb,amsxtra, mathrsfs,enumerate,relsize}
\usepackage{calrsfs}
\usepackage{array,booktabs,longtable}
\usepackage[us,12hr]{datetime}
\usepackage{hyperref}
 \hypersetup{colorlinks=true,citecolor=blue}
   \usepackage{graphicx}
     \usepackage[margin=1.5in]{geometry}
   \usepackage[font=scriptsize]{caption}
\usepackage{todonotes}
\usepackage{tikz}
\usepackage{tikz-cd}
\usetikzlibrary{positioning}
\CDat
\usepackage{scalerel}[2016/12/29]


\newtheorem{thm}{Theorem}[section]  
\newtheorem*{un-no-thm}{Theorem}
\newtheorem{cor}[thm]{Corollary}     
\newtheorem{lem}[thm]{Lemma}         
\newtheorem{prop}[thm]{Proposition}

\newtheorem{add}[thm]{Addendum}

\newtheorem{bigthm}{Theorem}

\newtheorem{bigcor}[bigthm]{Corollary}

\theoremstyle{definition}

\theoremstyle{definition}

\theoremstyle{definition}

\theoremstyle{remark}
\newtheorem{rem}[thm]{Remark}

\newtheorem*{out}{Outline}

\newtheorem{ex}[thm]{Example}

\DeclareMathOperator{\Top}{Top}
\DeclareMathOperator{\Sp}{Spectra}

\begin{document}
\title{On Poincaré spaces which are not manifolds}
\date{\today}
\author{John R.\ Klein}
\address{Wayne State University, Detroit, MI 48202}
\email{klein@math.wayne.edu}
\begin{abstract}
This expository note gives a short derivation of the classical Spivak
normal fibration formula for the total space of a spherical fibration
over a Poincar\'e duality space. We review the
associated smooth, PL, and topological manifold realization criteria
and illustrate them with examples from the stable stems. 
The principal application combines these criteria with the
completed solution of the Kervaire invariant problem to identify the
exceptional stems in a specified family of three cell Poincar\'e spaces.
\end{abstract}
\maketitle
\setlength{\parindent}{15pt}
\setlength{\parskip}{1pt plus 0pt minus 1pt}

\def\Sp{\text{\bf Sp}}
\def\vo{\varOmega}
\def\vs{\varSigma}
\def\smsh{\wedge}
\def\flush{\flushpar}
\def\id{\text{id}}
\def\dbslash{/\!\! /}
\def\codim{\text{\rm codim\,}}
\def\:{\colon}
\def\holim{\text{holim\,}}
\def\hocolim{\text{hocolim\,}}
\def\Bbb{\mathbb}
\def\bold{\mathbf}
\def\Aut{\text{\rm Aut}}
\def\cal{\mathcal}
\def\sec{\text{\rm sec}}
\def\gda{G\text{\rm -}\delta\text{\rm -}\alpha}
\def\PDD{\text{\rm pd\,}}
\def\PD{\text{\rm P}}
\def\stableto {\,\, \mapstochar \!\!\to}

\setcounter{tocdepth}{1}
\tableofcontents


\section{Introduction}

When speaking about Poincar\'e duality spaces, the author is frequently
asked whether there are examples that are not homotopy equivalent to
closed manifolds. This note grew out of that recurring question.

We give a short mapping cylinder derivation of the classical Spivak
fibration formula and use it to
organize explicit obstructions to manifold realization. The principal
application is Corollary~\ref{cor:kervaire-stems}, which combines the
classical realization criteria with the completed solution of the
Kervaire invariant problem.

The author claims little originality for this note and does not intend
to publish it in a peer-reviewed journal.

The examples illustrate the distinction between smooth and topological
manifold realization. They are total spaces of spherical fibrations with section.
Several are highly connected and have the integral cohomology ring
of a product of two spheres. Throughout, a {\it manifold realization} means
a homotopy equivalence to a closed manifold in the specified category.

The geometric input is the following spherical fibration analog of
the formula for the normal bundle of a sphere bundle with manifold base.
If $P$ is a Poincar\'e duality space, we
write $\nu_P$ for its Spivak normal fibration.

\begin{bigthm}[Spivak fibration formula {\cite[thm.~1]{Wissemann-Hartmann}}]\label{thm:main-formula}
Let $p\:E\to B$ be a spherical fibration $\xi$ with fiber
homotopy equivalent to $S^d$, where $d\ge2$ and $B$ is a finite
Poincar\'e space. Then there is a stable fiber homotopy equivalence
\[
p^*\nu_B\simeq p^*\xi+\nu_E.
\]
\end{bigthm}

Here, addition is fiberwise join of stable spherical fibrations.
We give a short derivation from Wall's normal fibration sum formula
\cite[chap.~11, p.~121]{Wall}.
If $s\:B\to E$ is a section, the formula implies
\[
\xi\simeq\nu_B-s^*\nu_E.
\]
Consequently, if $\nu_B$ admits a stable vector bundle reduction
but $\xi$ does not, then $E$ has no smooth manifold realization.
Replacing vector bundle reductions by $B\Top$-reductions gives
the corresponding topological obstruction. These statements are
Proposition~\ref{prop:reduction} and Addendum~\ref{add:topological-reduction} below.

\begin{rem} The author was aware of the Spivak fibration formula for more than twenty-five years, but first
became aware of Wissemann-Hartmann's paper \cite{Wissemann-Hartmann}
at a late stage in the preparation of this note.

Our proof first observes that the mapping cylinder inclusion $B\to M(p)$ is a 
Poincar\'e embedding. This enables us to relate the Spivak normal fibration of
$B$ to that of the Poincar\'e pair $(M(p),E)$.
Restricting to the boundary $E = \partial M(p)$ then yields the desired formula.

The proof in \cite[prop.~1.1 and cor.~2]{Wissemann-Hartmann} instead constructs
a Poincar\'e embedding $E\to B\times S^N$ and obtains the
mapping cylinder formula by doubling. An equivalent identity for
Poincar\'e pairs is used in \cite[\S3]{Klein_immersion}, and a
general product formula is proved in
\cite[thm.~I]{Klein_dualizing}.
\end{rem}

Table~\ref{tab:originality} summarizes the relation of the statements
below to earlier work. No priority is claimed for the precise
formulations of the corollaries and examples.
Proposition~\ref{prop:reduction} and
Addendum~\ref{add:topological-reduction} are formal consequences of
Theorem~\ref{thm:main-formula};
Remark~\ref{rem:compact-top-reductions} applies the results of
\cite{Klein-Williams_transfer}.

\begingroup
\small
\renewcommand{\arraystretch}{1.12}
\setlength{\LTpre}{\medskipamount}
\setlength{\LTpost}{\medskipamount}
\captionsetup{font=small}
\begin{longtable}{@{}>{\raggedright\arraybackslash}p{.24\linewidth}>{\raggedright\arraybackslash}p{\dimexpr.76\linewidth-2\tabcolsep\relax}@{}}
\caption{Scope and originality of the results}\label{tab:originality}\\
\toprule
\textbf{Result} & \textbf{Relation to earlier work}\\
\midrule
\endfirsthead
\caption[]{Scope and originality of the results (continued)}\\
\toprule
\textbf{Result} & \textbf{Relation to earlier work}\\
\midrule
\endhead
\midrule
\multicolumn{2}{r}{\footnotesize Continued on the next page}\\
\endfoot
\bottomrule
\endlastfoot
Theorem~\ref{thm:main-formula} &
The formula appears in \cite[thm.~1]{Wissemann-Hartmann}.
The mapping cylinder argument below is a short derivation from Wall's
normal fibration sum formula. No claim is made that this proof is new.\\
\addlinespace[.6em]
Theorem~\ref{thm:main-construction} &
The construction and realization criteria are classical
\cite{Smith_few_cells,Wissemann-Hartmann}.
The Kervaire invariant formulation uses the standard surgery
identification \cite[chap.~4]{Madsen-Milgram}.\\
\addlinespace[.6em]
Corollary~\ref{cor:kervaire-stems} &
A direct consequence of the classical criteria and the completed
Kervaire invariant problem \cite{Lin-Wang-Xu_Kervaire}.
This is the principal application of the note; no priority is claimed
for its precise formulation.\\
\addlinespace[.6em]
Theorem~\ref{thm:main-examples}(i) &
The classical five dimensional example
\cite[pp.~32--33]{Madsen-Milgram} and its fiberwise suspensions.
No new existence theorem is claimed.\\
\addlinespace[.6em]
Theorem~\ref{thm:main-examples}(ii)--(iii) &
Illustrations of the classical criteria using calculations in the stable stems:
the comparison involving $\kappa$ and the odd primary beta families.
These are presented as worked examples, without a claim of novelty.\\
\addlinespace[.6em]
Theorem~\ref{thm:main-examples}(iv) &
An application of the recently established class $\theta_6$
\cite[thm.~1.1]{Lin-Wang-Xu_Kervaire} to the classical construction.
It gives the $255$-dimensional example without introducing a new
obstruction or construction method.\\
\end{longtable}
\endgroup

\subsection{Examples}
In counting cells, we do not count the basepoint $0$-cell.
Write $\pi_r^S$ for the
$r$th stable homotopy group of spheres and let
$J\:\pi_r(O)\to\pi_r^S$ be the stable $J$-homomorphism.
For $r\equiv2\pmod4$, write
$K\:\pi_r^S\to\Bbb Z/2$ for the Kervaire invariant.

The fibration in the next theorem is obtained by the clutching construction. Let $F_d$
be the monoid of based self-homotopy equivalences of $S^d$.
If $r\ge1$, then by the Freudenthal suspension theorem we have
\[
\pi_r(F_{r+2})\cong\pi_{2r+2}(S^{r+2})
\xrightarrow[\cong]{\operatorname{st}}\pi_r^S.
\]
A representative $S^r\to F_{r+2}$ of the class corresponding to
$\alpha\in\pi_r^S$ is the clutching map. Since self-equivalences preserve
the fiber basepoint, the construction also gives a section; see
Section~\ref{subsec:clutching}.

\begin{bigthm}[Classical construction and realization criteria]\label{thm:main-construction}
For every $r\ge1$ and $\alpha\in\pi_r^S$, let
\[
S^{r+2}\longrightarrow E_\alpha\longrightarrow S^{r+1}
\]
be the spherical fibration obtained by clutching from $\alpha$.
It has a section, and its total space is an $r$-connected Poincar\'e
space of dimension
$2r+3$ with three cells and integral cohomology ring
$H^*(S^{r+1}\times S^{r+2};\Bbb Z)$.
If $\alpha\notin\operatorname{im}J$, then $E_\alpha$ has no smooth
manifold realization. 

For $r\ge6$ with
$r\equiv2\pmod4$, the space $E_\alpha$ has a topological manifold
realization if and only if $K(\alpha)=0$; when this holds, it has
a PL manifold realization.
\end{bigthm}

The construction and the smooth obstruction are given in
Corollary~\ref{cor:stable-examples} and the opening of
Section~\ref{sec:stable-stems}; the topological assertion is
Lemma~\ref{lem:kervaire-realization}.
Remark~\ref{rem:clutching-cell-model} describes the attaching map
of this three cell complex explicitly.
The underlying three cell construction and its smooth, PL, and
topological realization criteria, for sufficiently large fiber dimension,
are treated by Smith
\cite[thm., p.~279; props.~1--3]{Smith_few_cells}; see also
\cite[\S4.1.1 and rem.~9(a)]{Cornea_dynamics}.
The explicit choice of fiber dimension is convenient; the realization
criteria already apply in this range by
\cite[thms.~3 and~4]{Wissemann-Hartmann}.

The principal application of this note is the following consequence of
Theorem~\ref{thm:main-construction} and the completed solution of the
Kervaire invariant problem
\cite{Browder_Kervaire,Hill-Hopkins-Ravenel_Kervaire},
\cite[cor.~1.2]{Lin-Wang-Xu_Kervaire}.

\begin{bigcor}[Exceptional stems]\label{cor:kervaire-stems}
Suppose $r\ge6$ and $r\equiv2\pmod4$. There exists
$\alpha\in\pi_r^S$ such that $E_\alpha$ has no closed topological
manifold realization if and only if
\[
r\in\{6,14,30,62,126\}.
\]
For every other $r$ in this range, each $E_\alpha$ has a closed PL
manifold realization. If $\alpha\ne0$, it has no closed smooth
manifold realization.
\end{bigcor}

The examples in \cite[\S3]{Wissemann-Hartmann} concern manifold
homotopy types arising from fibrations that are not fiber homotopy
equivalent to sphere bundles. Here the examples emphasize failure
of manifold realization and the distinction between the smooth
and topological categories. Proposition~\ref{prop:reduction}
also records the immediate extension of the section obstruction
to Poincar\'e bases with reducible Spivak fibration.

For the Kervaire invariant one classes, $\kappa$, and the odd primary
beta elements,
the relevant group $\pi_r(O)$ vanishes. Their nonzero stable classes
therefore obstruct smooth realizations. The Kervaire invariant one
classes include $\nu^2$ and $\sigma^2$, whereas $K(\kappa)=0$ and the
Kervaire invariant vanishes on every class of odd order. These
calculations give the following collection of examples, recorded to
illustrate the realization criteria.

\begin{bigthm}[Illustrative examples]\label{thm:main-examples}
\begin{enumerate}[(i).]
\item For every $n\ge5$, there is a finite Poincar\'e space of
 dimension $n$ that is simply connected and has no topological
manifold realization.
\item There are two $14$-connected, $31$-dimensional Poincar\'e
spaces $X$ and $Y$, each with three cells and integral cohomology
ring
\[
H^*(X;\Bbb Z)\cong H^*(Y;\Bbb Z)\cong
H^*(S^{15}\times S^{16};\Bbb Z),
\]
such that neither has a smooth manifold realization, $X$ has no
topological manifold realization, and $Y$ has a PL manifold realization.
\item For every odd prime $p$, put $r_p=2p(p-1)-2$. There is an
$r_p$-connected closed PL manifold $M_p$ of dimension
$2r_p+3=4p(p-1)-1$ such that
\[
H^*(M_p;\Bbb Z)\cong
H^*(S^{r_p+1}\times S^{r_p+2};\Bbb Z)
\]
as graded rings, but $M_p$ is not homotopy equivalent to any closed
smooth manifold. The first three dimensions are $23$, $79$, and $167$.
\item There is a $126$-connected, $255$-dimensional Poincar\'e space
with three cells and integral cohomology ring
$H^*(S^{127}\times S^{128};\Bbb Z)$ which has no topological
manifold realization.
\end{enumerate}
\end{bigthm}

\begin{ex}[Madsen--Milgram]\label{ex:intro-madsen-milgram}
The classical example of \cite[pp.~32--33]{Madsen-Milgram} is a direct
application of Theorem~\ref{thm:main-formula}. It is given by a spherical
fibration $\xi$
\[
S^2\to E\to S^3
\]
with section $s$ and stable clutching class $0\ne\eta^2\in\pi_2^S$.
Since $\nu_{S^3}$ is stably trivial, Theorem~\ref{thm:main-formula} gives
$s^*\nu_E\simeq-\xi$. As $\pi_3(B\Top)=0$ \cite{Kirby-Siebenmann},
$\xi$ has no $B\Top$-reduction, and hence neither does $\nu_E$.
Thus $E$ is a simply connected, five dimensional Poincar\'e space
with no closed topological manifold realization.

Its twofold fiberwise suspension is the case $r=2$, $\alpha=\eta^2$
of Theorem~\ref{thm:main-construction}, up to fiber homotopy equivalence,
with fiber $S^4$ and total dimension $7$.
Iterated fiberwise suspensions preserve the obstruction and give
examples in every dimension at least five, proving
Theorem~\ref{thm:main-examples}(i). See Examples~\ref{ex:madsen-milgram}
and~\ref{ex:joins} for the construction and its suspensions.
\end{ex}

For Theorem~\ref{thm:main-examples}(ii), the spaces are $X=E_{\sigma^2}$ and $Y=E_\kappa$
from Example~\ref{ex:kappa}. Their different topological behavior is
detected by the Kervaire invariant. The class $\nu^2$ similarly gives
a $6$-connected, $15$-dimensional Poincar\'e space with the integral
cohomology ring of $S^7\times S^8$ and no topological manifold
realization (Example~\ref{ex:kervaire-classes}).
Part~(iii) uses the first beta element at each odd prime.
There are also infinite families at each fixed prime $p\ge5$:
at $p=5$, the class $\beta_t$ gives connectivity $48t-10$ and
dimension $96t-17$, for every $t\ge1$
(Example~\ref{ex:beta-family}).
Part~(iv) uses the Kervaire invariant one class in $\pi_{126}^S$
whose existence was proved by Lin--Wang--Xu
\cite[thm.~1.1]{Lin-Wang-Xu_Kervaire}.
Example~\ref{ex:kervaire-classes} also gives the corresponding
nonmanifold examples in dimensions $63$ and $127$.
\medskip

\noindent\emph{Use of a large language model.}
The author wrote the first draft of this paper without assistance from
a large language model. The author subsequently used OpenAI's Codex
to assist with revising the exposition,
locating and comparing references, examining mathematical arguments,
and developing examples and applications. The author takes full
responsibility for the mathematical content, references, and final
presentation of the paper.

\begin{out} Section~\ref{sec:definitions} gives the definitions and
conventions. Section~\ref{sec:spivak-formula} contains the proof of
Theorem~\ref{thm:main-formula}. Section~\ref{sec:reduction} establishes
the reduction criterion and applies it to the classical five dimensional
example and its fiberwise suspensions, proving
Theorem~\ref{thm:main-examples}(i). In Section~\ref{sec:stable-stems}
we complete the proof of Theorem~\ref{thm:main-construction}, prove
Corollary~\ref{cor:kervaire-stems}, and construct the highly connected
examples of Theorem~\ref{thm:main-examples}(ii)--(iv),
including the comparison in dimension $31$, the odd primary families,
and the Kervaire invariant one examples through dimension $255$.
\end{out}

\section{Definitions and conventions}\label{sec:definitions}

\subsection{Poincar\'e spaces and pairs}
An \emph{orientation sheaf} on a space $X$ is a local system
$\cal L$ of infinite cyclic groups. A connected space $X$ of
finite CW homotopy type is a \emph{Poincar\'e space of formal dimension
$n$} if it has an orientation sheaf and a fundamental class
$[X]\in H_n(X;\cal L )$ such that cap product gives isomorphisms
\[
H^i(X;\mathcal M)\xrightarrow{\;\cap[X]\;}
H_{n-i}(X;\cal L\otimes\mathcal M)
\]
for every locally constant sheaf $\mathcal M$ of abelian groups and every $i$.
The space is \emph{orientable} if $\cal L $ is constant.
The dimension always means the formal dimension.

A \emph{Poincar\'e pair} $(X,\partial X)$ of dimension $n$ has finite CW
homotopy type, with $X$ connected and $\partial X$ possibly empty. It has an
orientation sheaf $\cal L $ and a class
$[X]\in H_n(X,\partial X;\cal L )$ inducing
the duality isomorphisms
\begin{align*}
H^i(X;\mathcal M)&\xrightarrow{\;\cap[X]\;}
H_{n-i}(X,\partial X;\cal L \otimes\mathcal M),\\
H^i(X,\partial X;\mathcal M)&\xrightarrow{\;\cap[X]\;}
H_{n-i}(X;\cal L \otimes\mathcal M)
\end{align*}
for every local system $\mathcal M$ on $X$ and every $i$.
In addition, the boundary class $\partial[X]$ must give Poincar\'e
duality of dimension $n-1$ on each component of $\partial X$,
with orientation sheaf $\cal L|_{\partial X}$
\cite[def.~1.1]{Klein-Qin-Su}. A Poincar\'e
space without a specified boundary is understood to be closed.
Pairs are replaced by cofibrations when necessary.

\subsection{Spherical fibrations and fiberwise joins}
A \emph{$d$-spherical fibration} $\xi$ over $B$ is a fibration
$p\:E\to B$ whose fibers have the homotopy type of $S^d$.
Equivalence means fiber homotopy equivalence over $B$.
A \emph{section} is a map $s\:B\to E$ with $ps=\id_B$;
a homotopy section can be replaced by a section using the homotopy
lifting property.

For a $d$-spherical fibration with $d\ge1$ and connected base $B$
of finite CW homotopy type, the total space $E$ also has finite CW
homotopy type. In this setting, $E$ is a Poincar\'e space of dimension
$d+k$ if and only if $B$ is a Poincar\'e space of dimension $k$
\cite[cor.~F]{Klein_dualizing}.

If $\xi$ and $\eta$ have total spaces $E$ and $F$ over the same
base, their \emph{fiberwise join} $\xi+\eta$ is represented by the
double mapping cylinder of the projections
\[
E\longleftarrow E\times_B F\longrightarrow F,
\]
viewed over $B$ and replaced by a fibration if necessary.
Its fiber over $b$ is homotopy equivalent to $E_b*F_b$.
Thus the join of a $d$-spherical and an $e$-spherical fibration is
$(d+e+1)$-spherical. Joining $E\to B$ with the trivial $S^0$-fibration is
\emph{unreduced fiberwise suspension} $S_B E\to B$.

For $j\ge1$, write $S^j_B E$ for the total space of the join with
$B\times S^{j-1}\to B$, and set $S^0_B E=E$.
Thus $S^j_B E$ is the $j$-fold unreduced fiberwise suspension of
$E$ over $B$. Associativity of fiberwise join gives a fiber homotopy
equivalence
\[
S^i_B(S^j_B E)\simeq S^{i+j}_B E
\]
for $i,j\ge0$.
If $d\ge1$ and $B$ is a Poincar\'e space of dimension $k$, then
$S^j_B E$ is a Poincar\'e space of dimension $d+k+j$, and
$(M(p_j),S^j_B E)$ is a Poincar\'e pair of dimension $d+k+j+1$.
Here $p_j\:S^j_B E\to B$ is the projection and $M(p_j)$ is its
mapping cylinder.

Two spherical fibrations are \emph{stably equivalent} if suitable
iterated fiberwise suspensions are fiber homotopy equivalent.
Let $G$ denote the stable monoid of self-homotopy equivalences of
spheres. Stable spherical fibrations over $B$ are classified by
$[B,BG]$, an abelian group under fiberwise join. Trivial fibrations
represent zero, and $-\xi$ denotes the additive inverse of $\xi$.
Pullback preserves this group structure. All identities involving
Spivak normal fibrations in this paper are identities of stable classes.

\subsection{The Spivak normal fibration}
For a spherical fibration $\xi$ with fiber $S^{q-1}$, write
$S(\xi)$ for its total space and $D(\xi)$ for the mapping cylinder
of its projection. Its \emph{Thom space} is
$T(\xi)=D(\xi)/S(\xi)$.
If $X$ is an $n$-dimensional Poincar\'e space, a \emph{Spivak normal
fibration} consists of such a fibration $\nu_X$ over $X$, whose
orientation sheaf is that of $X$, together with a map
\[
c\:S^{n+q}\longrightarrow T(\nu_X)
\]
whose image of the fundamental class corresponds to $[X]$ under
the Thom isomorphism
\[
\widetilde H_{n+q}(T(\nu_X);\Bbb Z)
\cong H_n(X;\cal L ).
\]
Such data exist after stabilization, and the stable fiber homotopy
class of $\nu_X$ is unique; see \cite{Spivak} and
\cite[\S3]{Klein_Poincare}.

For a Poincar\'e pair $(X,\partial X)$, use the relative Thom space
$T(\nu)/T(\nu|_{\partial X})$ and the relative fundamental class
in the same definition. We denote the resulting stable fibration
over $X$ by $\nu_{(X,\partial X)}$. Its restriction to the boundary
is stably equivalent to $\nu_{\partial X}$
\cite[\S3]{Klein_immersion}.
For a smooth manifold, the Spivak fibration is the stable spherical
fibration of its normal bundle. In particular, $\nu_{S^k}$ is
stably trivial.

\subsection{Reductions and manifold realizations}
The maps
\[
BO\longrightarrow B\mathrm{PL}\longrightarrow B\Top
\longrightarrow BG
\]
forget successively the additional structure on stable bundles.
Here $B\mathrm{PL}$ and $B\Top$ classify stable PL and topological
Euclidean bundles. A \emph{stable vector bundle reduction} of
$\xi\:B\to BG$ is a homotopy lift through $BO\to BG$.
Reductions to $B\mathrm{PL}$ and $B\Top$ are defined in the same
way. The forgetful maps send Whitney sum to fiberwise join, so
the classes admitting a reduction form a subgroup of $[B,BG]$.
Reductions are also preserved under pullback.

A \emph{manifold realization} of a closed Poincar\'e space $X$
is a homotopy equivalence $M\simeq X$ with $M$ a closed manifold
in the indicated category. Its stable normal bundle supplies a
reduction of $\nu_X$ to $BO$, $B\mathrm{PL}$, or $B\Top$,
respectively. Consequently, failure of the corresponding reduction
obstructs a manifold realization. A reduction alone is not in
general sufficient; surgery obstructions may remain
\cite{Wall}.

\subsection{Clutching and the stable $J$-homomorphism}\label{subsec:clutching}
Let $F_d$ be the topological monoid of based self-homotopy
equivalences of $S^d$. Its classifying space $BF_d$ classifies
$d$-spherical fibrations with section. For $k\ge2$, after choosing
an identification of the fiber over the basepoint, the clutching
class of such a fibration over $S^k$ is
\[
\theta_\xi\in\pi_{k-1}(F_d)\cong\pi_{k+d-1}(S^d).
\]
A representative $S^{k-1}\to F_d$ specifies, in the fiber homotopy
sense, how to identify the trivial fibrations over the two hemispheres
of $S^k$. The maps preserve the fiber basepoint, so the basepoint
sections fit together to give a section of the resulting fibration.
Stabilization gives
\[
\theta_\xi^{\mathrm{st}}\in\pi_{k-1}^S,
\qquad
\pi_r^S=\operatorname*{colim}_{m}\pi_{r+m}(S^m).
\]
For $r\ge1$, the identification $\pi_{r+1}(BG)\cong\pi_r^S$
identifies the homomorphism induced by $BO\to BG$ with
\[
J\:\pi_r(O)\longrightarrow\pi_r^S.
\]
Thus the image of $J$ consists of the stable classes over
$S^{r+1}$ arising from vector bundles. We use
$K\:\pi_{4m+2}^S\to\Bbb Z/2$ for the Kervaire invariant of
the framed bordism class corresponding to a stable homotopy class
under the Pontryagin--Thom construction \cite{Browder_Kervaire}.

\begin{rem}[The attaching map]\label{rem:clutching-cell-model}
Let $\xi$ be a $d$-spherical fibration over $S^k$ with a chosen
section, where $k,d\ge2$. Under the clutching identification above,
represent $\theta_\xi$ by a map
$a_\xi\:S^{k+d-1}\to S^d$. With compatible choices of orientations,
the total space is homotopy equivalent to the three cell complex
\cite[\S3]{Wissemann-Hartmann}
\[
E\simeq(S^k\vee S^d)\cup_f D^{k+d},
\qquad
f=[\iota_B,\iota_F]+\iota_F\circ a_\xi.
\]
Here $\iota_B\:S^k\to S^k\vee S^d$ and
$\iota_F\:S^d\to S^k\vee S^d$ are the inclusions corresponding
to the section and the fiber. The Whitehead product is the attaching
map for the top cell of $S^k\times S^d$; the second term records
the clutching map.
For $E_\alpha$ in Theorem~\ref{thm:main-construction}, take
$k=r+1$ and $d=r+2$, with $a_\xi$ representing the unique homotopy
class stabilizing to $\alpha$. This gives the explicit model
\eqref{eq:stable-cell-model} below.
\end{rem}

\section{The Spivak fibration formula}\label{sec:spivak-formula}

Let $p\:E\to B$ be as in Theorem~\ref{thm:main-formula}, and put
$k=\dim B$. The total space $E$ is a Poincar\'e space of dimension
$d+k$ \cite[cor.~F]{Klein_dualizing}.

\begin{proof}[Proof of Theorem~\ref{thm:main-formula}]
Let $D=M(p)$ be the mapping cylinder of $p$, with boundary inclusion
$i\:E\to D$, inclusion $z\:B\to D$, and retraction
$r\:D\to B$. Thus $rz=\id_B$, $ri=p$, and $zp\simeq i$.
The fiberwise cone construction and the Thom isomorphism, together
with duality on $B$, make $(D,E)$ a Poincar\'e pair of dimension
$d+k+1$. Write $\nu_{(D,E)}$ for its Spivak normal fibration.

The inclusion $z$ is a codimension $d+1$ Poincar\'e embedding
with normal spherical fibration $\xi$ and with complement a collar on $E$.
To see this, split the mapping cylinder interval at its midpoint: the inner
half is a mapping cylinder neighborhood of $B$, and the outer half
is identified with $E\times I$, a collar of the boundary. After retracting this collar,
the associated homotopy pushout is
\[
\begin{tikzcd}
E\arrow[r,"\id"]\arrow[d,"p"'] & E\arrow[d,"i"]\\
B\arrow[r,"z"'] & D.
\end{tikzcd}
\]
The square uses the mapping cylinder homotopy $zp\simeq i$; the two
boundary faces of the collar are understood in the embedding data.
The normal fibration sum formula for Poincar\'e embeddings
\cite[chap.~11, p.~121]{Wall}, in its relative form for the ambient
Poincar\'e pair $(D,E)$, gives
\[
\nu_B\simeq\xi+z^*\nu_{(D,E)}.
\]
The relative form follows by applying the closed formula to the double
$D\cup_E D$ and restricting to either half.
Equivalently, $\nu_{(D,E)}\simeq r^*(\nu_B-\xi)$, the pair formula
recorded in \cite[cor.~2]{Wissemann-Hartmann}.
An equivalent formulation using Poincar\'e thickenings appears in
\cite[\S12, lem.~12.2 and cor.~12.3]{klein_surgery1}.
The same identity is used in
\cite[\S3, proof of thm.~3.1(2)]{Klein_immersion}, in the form
$\nu_P+\partial=\nu_V$: here $P$ is a Poincar\'e pair whose underlying
space $V$ is a closed Poincar\'e space, and $\partial$ is the spherical
fibration obtained from its boundary inclusion. In the present case
these are $(D,E)$, $B$, and $\xi$, respectively, under $D\simeq B$.
The boundary restriction of $\nu_{(D,E)}$ is stably $\nu_E$.
Pulling back along $p$ therefore gives
\begin{align*}
p^*\nu_B
&\simeq p^*\xi+(zp)^*\nu_{(D,E)}\\
&\simeq p^*\xi+i^*\nu_{(D,E)}\\
&\simeq p^*\xi+\nu_E.\qedhere
\end{align*}
\end{proof}

\begin{rem}
More generally, \cite[thm.~I and \S7]{Klein_dualizing} gives
\[
\nu_E\simeq\nu_p+p^*\nu_B
\]
for a fibration $F\to E\xrightarrow{p}B$ of connected finitely dominated
Poincar\'e spaces. Here $\nu_p$ is the fiberwise Spivak normal fibration
over $E$, whose restriction to $F$ is $\nu_F$. In the spherical case,
Theorem~\ref{thm:main-formula} identifies this term as
$\nu_p\simeq-p^*\xi$. Consequently, \cite[thm.~I]{Klein_dualizing} supplies
the general product formula, whereas the spherical fibration formula
above explicitly identifies its fiberwise term.
\end{rem}

The case of a trivial Spivak normal fibration on the base is an
immediate specialization.

\begin{cor}\label{cor:trivial-base}
If $\nu_B$ is trivializable, then there is a stable fiber homotopy
equivalence
\[
\nu_E\simeq -p^*\xi,
\]
where $-p^*\xi$ denotes the additive inverse of $p^*\xi$ in the group
of stable spherical fibrations over $E$.
\end{cor}

\section{Reduction criteria and examples}\label{sec:reduction}

Let $\xi$ be a $d$-spherical fibration over a finite Poincar\'e space
$B$, with projection $p\: E\to B$.
The following result will be used as a recognition principle for
deciding when the total space $E$ fails to admit a smooth or topological
manifold realization.

\begin{samepage}
\begin{prop}\label{prop:reduction}
Assume that
\begin{enumerate}[(i).]
\item $d\ge 2$,
\item $\xi$ has a section,
\item $\xi$ does not admit a reduction to a stable vector bundle, and
\item the Spivak normal fibration $\nu_B$ admits a reduction to a stable
vector bundle.
\end{enumerate}
Then $E$ does not have the homotopy type of a closed smooth manifold.
\end{prop}
\end{samepage}

\begin{proof}
Let $s\:B\to E$ be a section. Pulling the identity in
Theorem~\ref{thm:main-formula} back along $s$, and using
$s^*p^*=(p\circ s)^*=\id$, gives
\[
\nu_B\simeq\xi+s^*\nu_E,
\qquad\text{hence}\qquad
\xi\simeq\nu_B-s^*\nu_E.
\]
If $E$ had the homotopy type of a closed smooth manifold, its Spivak
normal fibration $\nu_E$ would admit a stable vector bundle reduction.
Hypothesis~(iv), together with preservation of reductions under pullback
and additive inverse, would then give a stable vector bundle reduction
of $\xi$, contradicting hypothesis~(iii).
\end{proof}

The same argument with $B\Top$-reductions proves the following.

\begin{add}\label{add:topological-reduction}
Assume hypotheses~(i) and~(ii) of Proposition~\ref{prop:reduction}.
If $\xi$ has no $B\Top$-reduction and $\nu_B$ has a $B\Top$-reduction,
then $E$ does not have the homotopy type of a closed topological manifold.
\end{add}

\begin{rem}\label{rem:categories}
Let $BO\to BG$ classify the stable spherical fibration associated to a
stable vector bundle. Hypothesis~(iii) says that there is no dotted lift
making the following diagram commute up to homotopy:
\[
\begin{tikzcd}
& BO\arrow[d]\\
B\arrow[ur,dashed]\arrow[r,"\xi"'] & BG.
\end{tikzcd}
\]
Hypothesis~(iv) holds in particular if $\nu_B$ is trivializable.
Failure of a $BO$-reduction alone does not imply the topological
conclusion of Addendum~\ref{add:topological-reduction}.
\end{rem}

\begin{ex}\label{ex:madsen-milgram}
With the clutching notation of Section~\ref{sec:definitions}, take
$k=3$, $d=2$, and let $\xi$ correspond to the nonzero element
\[
\theta_\xi=\eta_2\eta_3\in
\pi_2(F_2)\cong\pi_4(S^2)\cong\Bbb Z/2,
\]
where $\eta_2\:S^3\to S^2$ is the Hopf map and $\eta_3=\Sigma\eta_2$.
Its stabilization is the nonzero class
$\eta^2\in\pi_2^S\cong\Bbb Z/2$.
As $\pi_2(O)=0$, the five dimensional total space has no closed smooth
manifold realization. Its usual cellular model is
\[
E\simeq (S^2\vee S^3)
\cup_{[\iota_2,\iota_3]+\iota_2\circ\eta_2\eta_3}D^5,
\]
where $\iota_2$ and $\iota_3$ are the wedge inclusions. This is the
example described in \cite[p.~32]{Madsen-Milgram}; it is also reviewed
in \cite[intro. and \S4.6]{Klein_Poincare}, where it is identified
as the Gitler--Stasheff example.
For this example, the argument using the section to rule out a stable
vector bundle reduction was given by the author in the MathOverflow
answer \cite{Klein_MO_manifold}.

In fact $\pi_3(B\Top)=\pi_2(\Top)=0$. This follows from
$\Top/\mathrm{PL}\simeq K(\Bbb Z/2,3)$ and the $6$-connectivity of $\mathrm{PL}/O$
\cite{Kirby-Siebenmann}. Thus the nonzero stable class
$\xi\in\pi_3(BG)$ has no $B\Top$-reduction either.
Addendum~\ref{add:topological-reduction} therefore rules out a closed
topological manifold realization of $E$.

Twice suspending the fibers gives the clutching class
$\Sigma^2(\eta_2\eta_3)\in\pi_6(S^4)$, which corresponds to
$\eta^2$ under the Freudenthal isomorphism
$\pi_6(S^4)\cong\pi_2^S$.
Thus the twofold fiberwise suspension of this fibration is fiber
homotopy equivalent to $E_{\eta^2}\to S^3$, the case $r=2$ of
Theorem~\ref{thm:main-construction}.
\end{ex}

\begin{ex}\label{ex:joins}
Let $\xi$ be a $d$-spherical fibration over $B$ satisfying the
hypotheses of Proposition~\ref{prop:reduction}. For $j\ge0$, let
$\xi_j$ be its $j$-fold unreduced fiberwise suspension. Its total
space is $S^j_BE$, its fiber is $S^{d+j}$, and its total dimension
is $d+\dim B+j$.

Each $\xi_j$ has a section and the same stable class as $\xi$,
so it also satisfies the hypotheses of
Proposition~\ref{prop:reduction}. We obtain the family
\[
E,\quad S_B E,\quad S^2_B E,\quad\ldots
\]
with no closed smooth manifold realizations.
The same argument applies to Addendum~\ref{add:topological-reduction}.
Starting with $\xi$ from Example~\ref{ex:madsen-milgram}, we obtain simply connected
Poincar\'e spaces with no closed topological manifold realizations in
every dimension at least five.
\end{ex}

\begin{rem}\label{rem:compact-top-reductions}
Recall that a fibration  $F\to E\to B$ of homotopy finite spaces
admits  a {\it compact {\rm TOP} reduction} if there is
a fiber homotopy equivalence to a
fiber bundle of compact topological manifolds, whose fibers may have a boundary
\cite[defn.~1.1 and rem.~1.2]{Klein-Williams_transfer}. The fibration has a
{\it compact fiber smoothing} if it admits a fiber homotopy equivalence to
a smooth fiber bundle with compact manifold fibers, possibly with boundary.

The once fiberwise suspended fibration of Example~\ref{ex:madsen-milgram},
$S^3\to S_{S^3}E\to S^3$, is the example of
\cite[thm.~G]{Klein-Williams_transfer}. It admits a compact TOP reduction
but no compact fiber smoothing.
Every fibration in Theorem~\ref{thm:main-construction} also admits a
compact TOP reduction.

More precisely, the fibrations of Theorem~\ref{thm:main-construction}
admit a reduction of 
structure group to \[
\Top(S^{r+2}\times D^j,*),
\] 
for sufficiently large $j$, where $\ast$ is
the basepoint in the interior of $S^{r+2}\times D^j$
\cite[thms.~B and~E; rem.~11.1]{Klein-Williams_transfer}.

A compact TOP reduction does not imply a $B\Top$-reduction of the
stable spherical fibration or a closed topological manifold realization
of the total space.
\end{rem}

Let  $2\le k\le d-1$ and let
$\alpha\in\pi_{k-1}^S$ have nonzero image in the
cokernel of $J$. By the Freudenthal suspension theorem, the stabilization
homomorphism
\[
\operatorname{st}\:\pi_{k+d-1}(S^d)
\longrightarrow\pi_{k-1}^S
\]
is an isomorphism. Thus $\alpha$ determines a $d$-spherical fibration
$\xi_\alpha$ over $S^k$ with section.  
The following immediate corollary to Proposition \ref{prop:reduction}
is the nonmanifold implication of
\cite[cor.~3.1]{Wissemann-Hartmann}.

\begin{cor}\label{cor:stable-examples}
Assume $2\le k\le d-1$. If
$\alpha\in\pi_{k-1}^S$ has nonzero image in the
cokernel of
\[
J\:\pi_{k-1}(O)\longrightarrow\pi_{k-1}^S,
\]
then the total space of $\xi_\alpha$ does not have the homotopy type of
a closed smooth manifold.
\end{cor}

\section{Examples from stable homotopy groups}\label{sec:stable-stems}

We apply Corollary~\ref{cor:stable-examples} to elements of the stable stems,
including those supplied by the recently completed solution of the Kervaire
invariant problem. The additional topological assertions use the
Kervaire invariant and the category versions of the criterion in
\cite[thms.~3 and~4]{Wissemann-Hartmann}; they do not follow from
failure of vector bundle reduction alone.

Given $r\ge 1$ and
$\alpha\in\pi_r^S$, take $k=r+1$ and $d=r+2$, the smallest fiber
dimension allowed by Corollary~\ref{cor:stable-examples}. Let
\[
S^{r+2}\longrightarrow E_\alpha\longrightarrow S^{r+1}
\]
be the spherical fibration obtained by clutching from $\alpha$, as in
Theorem~\ref{thm:main-construction}. It has a section, and its total space is
an $r$-connected Poincar\'e space of dimension $2r+3$.
Choose a representative $a_\alpha\:S^{2r+2}\to S^{r+2}$ of the
unique homotopy class stabilizing to $\alpha$.
Remark~\ref{rem:clutching-cell-model} gives, with compatible choices
of orientations, the three cell complex
\begin{equation}\label{eq:stable-cell-model}
E_\alpha\simeq
(S^{r+1}\vee S^{r+2})
\cup_{[\iota_{r+1},\iota_{r+2}]+\iota_{r+2}\circ a_\alpha}
D^{2r+3}.
\end{equation}
Here the $\iota_j$ are the wedge inclusions.
The Whitehead product term gives the top dimensional cup product,
so
\[
H^*(E_\alpha;\Bbb Z)\cong
H^*(S^{r+1}\times S^{r+2};\Bbb Z)
\]
as graded rings. Nevertheless, $E_\alpha$ has no closed smooth
manifold realization whenever $\alpha\notin\operatorname{im}J$.

For the examples below it is also useful to record the topological
criterion explicitly. Let
$K\:\pi_{4m+2}^S\to\Bbb Z/2$ denote the Kervaire invariant.

\begin{lem}\label{lem:kervaire-realization}
Suppose $r\ge 6$ and $r\equiv 2\pmod 4$. Then $E_\alpha$ has the
homotopy type of a closed topological manifold if and only if
$K(\alpha)=0$. In that case it also has the homotopy type of a
closed PL manifold.
\end{lem}

\begin{proof}
For $\mathrm{CAT}=\mathrm{PL}$ or $\Top$, the fibration
$G/\mathrm{CAT}\to B\mathrm{CAT}\to BG$ gives the exact sequence
\[
\pi_r(\mathrm{CAT})\longrightarrow\pi_r(G)\cong\pi_r^S
\longrightarrow\pi_r(G/\mathrm{CAT})\cong\Bbb Z/2.
\]
The last map is the framed surgery obstruction, hence the Kervaire
invariant. The surgery identification and its PL and topological
versions are reviewed in \cite[cor.~2.24, rem.~2.25, and chap.~4]{Madsen-Milgram}.
Thus the stable class of the spherical fibration over $S^{r+1}$
admits a $B\mathrm{CAT}$-reduction exactly when $K(\alpha)=0$.
The result follows from \cite[thms.~3 and~4]{Wissemann-Hartmann}.
\end{proof}

\begin{proof}[Proof of Corollary~\ref{cor:kervaire-stems}]
The completed solution of the Kervaire invariant problem states that
$K\:\pi_r^S\to\Bbb Z/2$ is nonzero precisely for
$r\in\{2,6,14,30,62,126\}$, among stems congruent to $2$ modulo $4$.
The nonexistence results are due to Browder
\cite{Browder_Kervaire} and Hill--Hopkins--Ravenel
\cite{Hill-Hopkins-Ravenel_Kervaire}; Lin--Wang--Xu settled the last
existence case, $r=126$. The full list, including the earlier existence
results, is recorded in \cite[cor.~1.2]{Lin-Wang-Xu_Kervaire}.
For $r\ge6$, Lemma~\ref{lem:kervaire-realization} identifies failure
of topological manifold realization with $K(\alpha)=1$.
Such an $\alpha$ exists exactly in the five stems listed in the
statement. In every other stem in the indicated range, $K$ vanishes
identically, and the same lemma gives a PL manifold realization for
each $E_\alpha$.
Finally, Bott periodicity gives $\pi_r(O)=0$ for $r\equiv2\pmod4$,
so any nonzero $\alpha$ is outside $\operatorname{im}J$.
Corollary~\ref{cor:stable-examples} gives the smooth obstruction.
\end{proof}

\begin{ex}[Kervaire invariant one classes]\label{ex:kervaire-classes}
For $2\le j\le6$, put $r_j=2^{j+1}-2$ and choose
$\theta_j\in\pi_{r_j}^S$ with $K(\theta_j)=1$.
Such classes exist by \cite[cor.~1.2]{Lin-Wang-Xu_Kervaire}.
For the first two, we may take $\theta_2=\nu^2$ and
$\theta_3=\sigma^2$; see \cite{Browder_Kervaire} and
\cite[thm.~3.4.4 and tab.~A3.3]{Ravenel_green}.
The space $E_{\theta_j}$ is $r_j$-connected, has dimension $2r_j+3$,
and has three cells with integral cohomology ring
\[
H^*(E_{\theta_j};\Bbb Z)\cong
H^*(S^{r_j+1}\times S^{r_j+2};\Bbb Z).
\]
By Lemma~\ref{lem:kervaire-realization}, none of these spaces has a
closed topological manifold realization. The resulting examples are
\[
\begin{array}{c|c|c}
\text{class} & \text{connectivity} & \text{dimension}\\ \hline
\theta_2=\nu^2 & 6 & 15\\
\theta_3=\sigma^2 & 14 & 31\\
\theta_4 & 30 & 63\\
\theta_5 & 62 & 127\\
\theta_6 & 126 & 255.
\end{array}
\]
In particular, the class whose existence is established in
\cite[thm.~1.1]{Lin-Wang-Xu_Kervaire} gives a fibration with section
\[
S^{128}\longrightarrow E_{\theta_6}\longrightarrow S^{127}.
\]
Its total space is a $126$-connected, $255$-dimensional Poincar\'e
space with the integral cohomology ring of $S^{127}\times S^{128}$
and no closed topological manifold realization. This proves
Theorem~\ref{thm:main-examples}(iv).
These examples increase the connectivity, whereas fiberwise
suspension of Example~\ref{ex:madsen-milgram} retains the base $S^3$.
\end{ex}

\begin{ex}[Two different outcomes in dimension $31$]\label{ex:kappa}
The full stable stem $\pi_{14}^S$ is entirely $2$-primary,
with decomposition \cite[chaps.~X and~XIII]{Toda}
\[
\pi_{14}^S\cong
(\Bbb Z/2)\{\sigma^2\}\oplus(\Bbb Z/2)\{\kappa\},
\qquad K(\sigma^2)=1,\quad K(\kappa)=0.
\]
The values of the Kervaire invariant follow from Browder's detection
theorem \cite{Browder_Kervaire} and the computation recorded
in \cite[tab.~A3.3]{Ravenel_green}.
Equivalently, $\kappa$ is the unique nonzero element in the kernel
of $K\:\pi_{14}^S\to\Bbb Z/2$.

Both $E_{\sigma^2}$ and $E_\kappa$ are $14$-connected,
$31$-dimensional Poincar\'e spaces with three cells and integral
cohomology ring $H^*(S^{15}\times S^{16};\Bbb Z)$.
Since $\pi_{14}(O)=0$, neither has the homotopy type of a closed
smooth manifold. However, Lemma~\ref{lem:kervaire-realization}
shows that $E_\kappa$ has a closed PL manifold realization,
whereas $E_{\sigma^2}$ has no closed topological manifold realization.
Thus even within this small collection of three cell complexes,
the cohomology ring and connectivity do not distinguish the two
realization problems.
\end{ex}

\begin{ex}[Odd primary families]\label{ex:beta-family}
For every odd prime $p$, the first beta element is a nonzero
class of order $p$,
\[
\beta_1\in\pi_{r_p}^S,
\qquad r_p=2p(p-1)-2;
\]
see \cite[thm.~1.1.14]{Ravenel_green}.
Since $r_p\equiv 2\pmod 4$, Bott periodicity gives
$\pi_{r_p}(O)=0$. Corollary~\ref{cor:stable-examples} therefore gives
an $r_p$-connected Poincar\'e space $E_{\beta_1}$ of dimension
\[
2r_p+3=4p(p-1)-1
\]
with the integral cohomology ring of
$S^{2p(p-1)-1}\times S^{2p(p-1)}$ and no closed smooth manifold
realization. On the other hand, $K(\beta_1)=0$, since a homomorphism
to $\Bbb Z/2$ kills every element of odd order.
Lemma~\ref{lem:kervaire-realization} gives a closed PL manifold
realization in every case. The first three cases are
\[
\begin{array}{c|c|c|c}
p & \text{base} & \text{fiber} & \dim E_{\beta_1}\\ \hline
3 & S^{11} & S^{12} & 23\\
5 & S^{39} & S^{40} & 79\\
7 & S^{83} & S^{84} & 167.
\end{array}
\]

There are also infinite families at each fixed prime $p\ge 5$.
Smith's beta family theorem gives a nonzero class $\beta_t$ of
order $p$ in degree
\[
r_t=2(p^2-1)t-2p,\qquad t\ge 1
\]
\cite[thm.~1.3.15]{Ravenel_green}. Again $r_t\equiv2\pmod4$,
so the same argument gives closed PL manifold homotopy types of
dimension $2r_t+3$ and connectivity $r_t$ with no closed smooth
manifold realizations. For example, at $p=5$ these dimensions
are $96t-17$ and the connectivities are $48t-10$.
The obstruction in these examples is odd primary, while the
topological obstruction in Example~\ref{ex:kervaire-classes} is the
Kervaire invariant modulo two.
\end{ex}

\enlargethispage{3\baselineskip}

\end{document}